# Mean Curvature Flows of Lagrangian Submanifolds with Convex Potentials

Knut Smoczyk and Mu-Tao Wang

October 10, 2002

**Abstract**

This article studies the mean curvature flow of Lagrangian submanifolds. In particular, we prove the following global existence and convergence theorem: if the potential function of a Lagrangian graph in $T^{2n}$ is convex, then the flow exists for all time and converges smoothly to a flat Lagrangian submanifold. We also discuss various conditions on the potential function that guarantee global existence and convergence.

## 1 Introduction

The mean curvature flow is an evolution process under which a submanifold evolves in the direction of its mean curvature vector. It can be considered as the gradient flow of the area functional in the space of submanifolds. The critical points of the area functional are minimal submanifolds.

In mirror symmetry, a distinguished class of minimal submanifolds called "special Lagrangians" are desirable in any complex $n$ dimensional Calabi-Yau manifold with a parallel $(n, 0)$ form $\Omega$. A special Lagrangian is calibrated by $\operatorname{Re}\Omega$, which means $*\operatorname{Re}\Omega = 1$, where $*$ is the Hodge $*$ operator on the submanifold. A simple derivation using Stokes' Theorem shows a special Lagrangian minimizes area in its homology class. To produce special Lagrangians, it is thus natural to consider the mean curvature flow. We remark that the existence of Lagrangian minimizers in Kähler-Einstein surfaces are proved by Schoen-Wolfson[8] using variational method.

It is conjectured by Thomas and Yau in [16] that a stable Lagrangian isotopy class in a Calabi-Yau manifold contains a smooth special Lagrangian and the deformation process can be realized by the mean curvature flow. One of the stability condition is in terms of the range of $*\text{Re}\,\Omega$. In [18] (see also [22]), the second author proves the following regularity theorem,

**Theorem 1.1** *Let $(X, \Omega)$ be a Calabi-Yau manifold and $\Sigma$ be a compact Lagrangian submanifold. If $*Re\,\Omega > 0$ on $\Sigma$, the mean curvature flow of $\Sigma$ does not develop any type I singularity.*

In particular, this theorem implies no neckpinching will occur in the flow. We remark that without this condition, neckpinching is possible by an example of Schoen-Wolfson [9]. It is thus of great interest to identify initial conditions that guarantee the long time existence and convergence of the flow.

The mean curvature flow of Lagrangian surfaces in four-manifolds are studied in [15] and [19]. The first author [15] proves the long-time existence and smooth convergence theorem for graphs of area preserving diffeomorphisms in the non-positive curvature case assuming an angle condition. In [19], the second author proves the long-time existence for graphs of area preserving diffeomorphisms between Riemann surfaces and uniform convergence when the diffeomorphism is homotopic to identity (smooth convergence for spheres). This gives a natural deformation retract of the group of symplectomorphism of Riemann surfaces. The maximum principle for parabolic equations is important in both papers [19] and [15]. The new ingredient in [19] is the blow-up analysis of the mean curvature flow developed in [18]. This has been applied to prove long-time existence and convergence theorems for general graphic mean curvature flows in arbitrary dimension and codimension in [20].

In this article, we prove the following global existence and convergence theorem in arbitrary dimension.

**Theorem A** *Let $\Sigma$ be a Lagrangian submanifold in $T^{2n}$. Suppose $\Sigma$ is the graph of $f : T^n \mapsto T^n$ and the potential function $u$ of $f$ is convex. Then the mean curvature flow of $\Sigma$ exists for all time and converges smoothly to a flat Lagrangian submanifold.*

$u$ is only a locally defined function and the convexity of $u$ will be defined more explicitly in §2. The mean curvature flow can be written locally as a fully non-linear parabolic equation for the potential $u$.

$$\frac{du}{dt} = \frac{1}{\sqrt{-1}} \ln \frac{\det(I + \sqrt{-1}D^2 u)}{\sqrt{\det(I + (D^2 u)^2)}} \quad (1.1)$$

The Dirichlet problem for the elliptic equation of (1.1) was solved by Caffarelli, Nirenberg and Spruck in [1].

The general existence theorem in [20] specialized to the Lagrangian case holds under the assumption that $\prod(1 + \lambda_i^2) < 4$ where $\lambda_i's$ are eigenvalues of $D^2 u$. The first author proves the convexity of $u$ is preserved in [14] and he also shows the existence and convergence theorem assuming $u$ is convex and the eigenvalues of $D^2 u$ are less than one. The method in [14] indeed implies stronger results.

The core of the proof is to get control of $D^2 u$. It is interesting that there are two ways to interpret $D^2 u$. First we can identify it with a symmetric two-tensor on the submanifold $\Sigma$. One can then calculate the evolution equation with respect to the rough Laplacian on symmetric two tensors. Applying Hamilton's maximum principle [3] shows that the subset of positive definite symmetric two-tensors is preserved along the flow. A stronger positivity gives the uniform $C^2$ bound of $u$.

On the other hand, since $f = \nabla u$, the graph of $D^2 u = df$ is the tangent space of the graph of $f$. Recall the Gauss map for a submanifold assigns each point to its tangent space. It was proved in [23] that the Gauss map of any mean curvature flow is a harmonic map heat flow and thus any convex region of the Grassmannian is preserved along the flow. From this we identify an expression of $D^2 u$ that satisfies the maximum principle and this also gives uniform $C^2$ bound of $u$. The convergence part uses Krylov's $C^{2,\alpha}$ estimate [6] for nonlinear parabolic equations.

Since the geometry of a Lagrangian submanifold is invariant under the unitary group $U(n)$, this gives other equivalent conditions to the convexity of $u$ that also imply global existence and convergence. This is explained in §5.

The first author would like to thank J. Jost, S.-T. Yau, G. Huisken and K. Ecker for many helpful discussions and suggestions. The second author would like to thank D. H. Phong and S.-T. Yau for their constant encouragement and support. He has benefitted greatly from conversations with B. Andrews, T. Ilmanen, R. Hamilton and J. Wolfson.

The second author proved Theorem A in an earlier preprint of the same title. The first author then pointed out that the first part of Theorem 1.3 in [14] applied to $S(V,W) = -\langle JV, \bar{W}\rangle$ says that convexity is preserved and that $u$ is uniformly bounded in $C^2$ because all induced metrics are uniformly equivalent. In this case, condition (1.7) in the second part of that theorem becomes redundant since it was only needed to imply convexity and one obtains long time existence and convergence giving a different proof of Theorem A. However, this fact was not mentioned explicitly in the paper. The $C^2$ estimate in the earlier preprint of the second author uses the geometry of the Lagrangian Grassmannian and this also leads to the observation of the equivalence conditions under the $U(n)$ action in §5. The authors thus decided to write a joint paper to incorporate both approaches and clarify various assumptions.

## 2 Preliminary

We first derive the evolution equation of $f = \nabla u$ from the equation of the potential function $u$. For more material on the special Lagrangian equation, we refer to Harvey-Lawson [4].

**Definition 2.1** *Let $\Omega \subset \mathbb{R}^n$ be a domain. $u : \Omega \times [0,T) \mapsto \mathbb{R}$ is said to satisfy the special Lagrangian evolution equation if*

$$\frac{du}{dt} = \frac{1}{\sqrt{-1}} \ln \frac{\det(I + \sqrt{-1}D^2 u)}{\sqrt{\det(I + (D^2 u)^2)}} \tag{2.1}$$

It is not hard to check $\frac{\det(I+\sqrt{-1}D^2 u)}{\sqrt{\det(I+(D^2 u)^2)}}$ is a unit complex number, so the right hand side is always real.

**Proposition 2.1** *Let $u_i = \frac{\partial u}{\partial x^i}$, then $u_i$ satisfies the following evolution equation.*

$$\frac{du_i}{dt} = g^{jk} u_{ijk} \tag{2.2}$$

*where $g^{jk} = g_{jk}^{-1}$ and $g_{jk} = \delta_{jk} + u_{jl}u_{kl}$*

*Proof.* Use the formula $(\ln \det A)' = A'_{ij} A^{ji}$, we compute

$$\frac{du_i}{dt} = u_{ijk}(I + \sqrt{-1}D^2 u)^{-1}_{kj} - \frac{1}{2\sqrt{-1}}(u_{lji}u_{lk} + u_{lj}u_{lki})(I + (D^2 u)^2)^{-1}_{kj}$$

It is not hard to check that

$$(I + \sqrt{-1}D^2 u)^{-1} = (I + (D^2 u)^2)^{-1} - \sqrt{-1}D^2 u(I + (D^2 u)^2)^{-1}$$

Therefore

$$\frac{du_i}{dt} = u_{ijk}(I+(D^2u)^2)^{-1}_{kj} - \sqrt{-1}u_{ijk}u_{kl}(I+(D^2u)^2)^{-1}_{lj} + \sqrt{-1}u_{lji}u_{lk}(I+(D^2u)^2)^{-1}_{kj}$$

where we use symmetry in $k, j$ in the last equality, the last two terms cancel and the proposition is proved. $\square$

The right hand side in (2.2) is the mean curvature form $H_i = g^{jk}h_{ijk}$, i.e. the trace of the second fundamental tensor $h_{ijk}$ because in our local coordinates we have $h_{ijk} = u_{ijk}$. It is well known that the mean curvature form $H$ is closed. Locally (e.g. see section 2.6 in [13]) $H$ can be expressed by the differential $d\alpha$ of the Lagrangian angle $\alpha = \frac{1}{\sqrt{-1}} \ln \frac{\det(I+\sqrt{-1}D^2 u)}{\sqrt{\det(I+(D^2u)^2)}}$, i.e. the right hand side in (2.1). Then we can give another proof of (2.2) by

$$\frac{d}{dt}du = d\frac{du}{dt} = d\alpha = H.$$

Equation (2.2) is indeed the nonparametric form of a graphic mean curvature flow, see [24] or [20] for the derivation of the general case. The graph of $\nabla u$ is then a Lagrangian submanifold in $\mathbb{C}^n \cong \mathbb{R}^n \oplus \mathbb{R}^n$ evolving by the mean curvature flow. It is well-known that being Lagrangian is preserved along the mean curvature flow, see for example [11] or [13]. The complex structure $J$ on $\mathbb{C}^n$ is chosen so that the second summand $\mathbb{R}^n$ is the image under $J$ of the first summand. The equation (2.2) is equivalent up to tangential diffeomorphisms to the original flow.

Now suppose $f : T^n \times [0,T) \mapsto T^n$ is given so that the graph of $f$ is a Lagrangian submanifold moved by mean curvature flow in $T^{2n} \cong T^n \times T^n$. The tangent space of $T^{2n}$ is identified with $\mathbb{C}^n \cong \mathbb{R}^n \oplus \mathbb{R}^n$. Now the differential $df$ is a linear map from the first $\mathbb{R}^n$ to the second $\mathbb{R}^n$, so is the complex structure $J$. The Lagrangian condition implies the bilinear form $\langle df(\cdot), J(\cdot) \rangle$ is symmetric. This implies there is a locally defined potential function $u$ of $f$. We shall identify $D^2 u$ with the bilinear form $\langle df(\cdot), J(\cdot) \rangle$.

**Definition 2.2** *The eigenvalues of $D^2 u$ are the eigenvalues of the symmetric bilinear form $\langle df(\cdot), J(\cdot)\rangle$. $u$ is convex if $\langle df(v), J(v)\rangle > 0$ for any $v \in \mathbb{R}^n$.*

Therefore an eigenvalue $\lambda$ of $D^2 u$ satisfies $df(v) = \lambda J(v)$ for some nonzero $v \in \mathbb{R}^n$. It is not hard to check that by integration the potential $u$ satisfies the special Lagrangian evolution equation locally.

## 3 Lagrangian Grassmannians

Let $LG(n)$ denote the Lagrangian Grassmannian of all Lagrangian subspaces of $\mathbb{C}^n$. Let $O$ be a base point in the Lagrangian Grassmannian $LG(n)$. $O$ represents a Lagrangian subspace in $\mathbb{C}^n$. A local coordinate chart near $O$ is parametrized by $\mathfrak{S}$, the space of $n \times n$ symmetric matrices of the form $Z = [z_{ij}]_{i,j=1\cdots n}$. They represent the collection of all Lagrangian subspaces that can be written as a graph over $O$. By [26], the invariant metric on $LG(n)$ is given by
$$ds^2 = Tr[(I+Z^2)^{-1} dZ]^2$$

Let $\phi(z_{ij})$ be a function on $\mathfrak{S}$. Given any geodesic $z_{ij}(s)$ let $\mathfrak{p}(s) = \phi(z_{ij}(s))$. Then $\phi$ is convex if and only if $\ddot{\mathfrak{p}}(s) \geq 0$.

We shall check convexity at an arbitrary point $P$ in $\mathfrak{S}$. $P$ is spanned by

$$\{a_i + z_{ij}(0) J a_j\}_{i=1\cdots n} \tag{3.1}$$

where $\{a_i\}_{i=1\cdots n}$ and $\{J a_i\}_{i\cdots n}$ form orthonormal bases for $O$ and $O^\perp \cong JO$, respectively. We may assume the bases are chosen so that $z_{ij}(0) = \lambda_i \delta_{ij}$ by diagonalization. By [26], a geodesic through $P$ parametrized by arc length is given as $P_s$ spanned by $\{a_i + z_{ij}(s) J a_j\}_{i=1\cdots n}$ such that $Z = [z_{ij}(s)]$ is a $n \times n$ matrix which satisfies the following ordinary differential equation.

$$\ddot{Z} - 2\dot{Z} Z (I + Z^2)^{-1} \dot{Z} = 0 \tag{3.2}$$

Note that $Z = Z^T$ for Lagrangian Grassmanians.

In the rest of the section, we develop a criterion to check when a function defined in terms of the eigenvalues of $z_{ij}$ is a convex function on $LG(n)$. Recall $\mathfrak{p}(s) = \phi(z_{ij}(s))$, so we have:

$$\ddot{\mathfrak{p}} = \frac{\partial^2 \phi}{\partial z_{ij} \partial z_{mn}} \dot{z}_{ij} \dot{z}_{mn} + \frac{\partial \phi}{\partial z_{ij}} \ddot{z}_{ij} \tag{3.3}$$

We assume $\phi$ is given by the eigenvalues $\lambda_k$ of $z_{ij}$. The second author learned the following formula from Ben Andrews.

**Lemma 3.1** *Given any function $\phi$ of $\lambda_k$, then*

$$\frac{\partial \phi}{\partial z_{ij}} = \frac{\partial \phi}{\partial \lambda_i} \delta_{ij} \tag{3.4}$$

$$\frac{\partial^2 \phi}{\partial z_{ij} \partial z_{mn}} \dot{z}_{ij} \dot{z}_{mn} = \sum_{k,l} \frac{\partial^2 \phi}{\partial \lambda_k \partial \lambda_l} \dot{z}_{kk} \dot{z}_{ll} + \sum_{k \neq l} \frac{\frac{\partial \phi}{\partial \lambda_k} - \frac{\partial \phi}{\partial \lambda_l}}{\lambda_k - \lambda_l} \dot{z}_{kl}^2 \tag{3.5}$$

Plug these equations into equation (3.3), we obtain

$$\ddot{\mathfrak{p}} = \sum_{k,l} \frac{\partial^2 \phi}{\partial \lambda_k \partial \lambda_l} \dot{z}_{kk} \dot{z}_{ll} + \sum_{k \neq l} \frac{\frac{\partial \phi}{\partial \lambda_k} - \frac{\partial \phi}{\partial \lambda_l}}{\lambda_k - \lambda_l} \dot{z}_{kl}^2 + \sum_i \frac{\partial \phi}{\partial \lambda_i} \ddot{z}_{ii} \tag{3.6}$$

Replace $\ddot{z}_{ij}$ by the geodesic equation (3.2),

$$\ddot{\mathfrak{p}} = \sum_{k,l} \frac{\partial^2 \phi}{\partial \lambda_k \partial \lambda_l} \dot{z}_{kk} \dot{z}_{ll} + \sum_{k \neq l} \frac{\frac{\partial \phi}{\partial \lambda_k} - \frac{\partial \phi}{\partial \lambda_l}}{\lambda_k - \lambda_l} \dot{z}_{kl}^2 + \sum_{i,p,k,l} \frac{\partial \phi}{\partial \lambda_i} 2\dot{z}_{ip} z_{pk} (\delta_{kl} + z_{km} z_{lm})^{-1} \dot{z}_{li}$$

Taking into account that $z_{ij}(0) = \lambda_i \delta_{ij}$, we obtain

**Proposition 3.1** *Let $\phi = \phi(\lambda_i)$ be a function on $\mathfrak{S}$. Along any geodesic $z_{ij}(s)$ on the Lagrangian Grassmannian with $z_{ij}(0) = \lambda_i \delta_{ij}$, at $s = 0$ we have*

$$\ddot{\mathfrak{p}} = \sum_{k,l} \frac{\partial^2 \phi}{\partial \lambda_k \partial \lambda_l} \dot{z}_{kk} \dot{z}_{ll} + \sum_{k \neq l} \frac{\frac{\partial \phi}{\partial \lambda_k} - \frac{\partial \phi}{\partial \lambda_l}}{\lambda_k - \lambda_l} (\dot{z}_{kl})^2 + \sum_{k,l} \frac{\partial \phi}{\partial \lambda_l} \left(\frac{2\lambda_k}{1+\lambda_k^2}\right) (\dot{z}_{kl})^2 \tag{3.7}$$

We apply this formula to the following function which corresponds to $-\ln \sqrt{\det(I + (D^2 u)^2)}$.

$$\phi_0 = -\frac{1}{2} \ln \prod (1 + \lambda_i^2)$$

We compute

$$\frac{\partial \phi_0}{\partial \lambda_k} = \frac{-\lambda_k}{1+\lambda_k^2}$$

Thus
$$\frac{\frac{\partial \phi_0}{\partial \lambda_k} - \frac{\partial \phi_0}{\partial \lambda_l}}{\lambda_k - \lambda_l} = \frac{\lambda_k \lambda_l - 1}{(1+\lambda_k^2)(1+\lambda_l^2)}$$

for $k \neq l$.
and
$$\frac{\partial^2 \phi_0}{\partial \lambda_k \partial \lambda_l} = \frac{-1+\lambda_k^2}{(1+\lambda_k^2)^2}\delta_{kl}$$

Therefore
$$\ddot{\mathfrak{p}} = \sum_k \frac{-1+\lambda_k^2}{(1+\lambda_k^2)^2}(\dot{z}_{kk})^2 + \sum_{k \neq l} \frac{\lambda_k\lambda_l - 1}{(1+\lambda_k^2)(1+\lambda_l^2)}(\dot{z}_{kl})^2 + \sum_{k,l} \frac{-2\lambda_k\lambda_l}{(1+\lambda_k^2)(1+\lambda_l^2)}(\dot{z}_{kl})^2$$

We obtain the general formula by splitting the last term according to $k = l$ or $k \neq l$.

$$\ddot{\mathfrak{p}} = \sum_k \frac{-1}{(1+\lambda_k^2)}(\dot{z}_{kk})^2 - \sum_{k \neq l} \frac{\lambda_k\lambda_l + 1}{(1+\lambda_k^2)(1+\lambda_l^2)}(\dot{z}_{kl})^2 \qquad (3.8)$$

So if $\lambda_k\lambda_l > -1$ for $k \neq l$, we have concavity.

**Proposition 3.2** *The function $\phi_0 = -\frac{1}{2}\ln\prod(1+\lambda_i^2)$ is a concave function on the subset of symmetric matrices where $\lambda_k\lambda_l > -1$ for $k \neq l$ of the Lagrangian Grassmannian.*

# 4 Proof of Theorem A

First we prove the convexity of $u$ is preserved as long as the flow exists smoothly. This is also proved in [14]. Consider the parametric version of the Lagrangian mean curvature flow $F : \Sigma \times [0, T) \mapsto T^{2n}$. The tangent space of $T^{2n}$ is identified with $\mathbb{C}^n \cong \mathbb{R}^n \oplus \mathbb{R}^n$ and the complex structure $J$ maps to first real space to the second summand. Let $\pi_1$ and $\pi_2$ denote the projection onto the first and second summand in the splitting. Define

$$S(X, Y) = \langle J\pi_1(X), \pi_2(Y) \rangle$$

for any $X, Y \in \mathbb{C}^n \cong T(T^{2n})$ as a two-tensor on $T^{2n}$. $S(X, Y)$ is a symmetric for any $X, Y$ in a Lagrangian subspace of $\mathbb{C}^n$. This is because $\omega(X, Y) = \langle J(\pi_1(X) + \pi_2(X)), \pi_1(Y) + \pi_2(Y) \rangle = 0$.

When $F$ is given as the graph of $f = T^n \mapsto T^n$, i.e. $F = (x, f(x))$. $df$ is a linear map, $df : \mathbb{R}^n \mapsto \mathbb{R}^n$. We have $dF(v) = v + df(v)$, $\pi_1(dF(v)) = v$ and $\pi_2(dF(v)) = df(v)$. Therefore $F^*S$ becomes a symmetric two-tensor and is the same as $\langle Jv, df(v) \rangle = \langle df(v), Jv \rangle$. $f$ has a locally defined potential $u$. By definition 2.2, the positive definiteness of $S$ is the same as the convexity of $u$.

Now we recall the general evolution equation for the pull back of a parallel two-tensor of the ambient space from [18] (§2, equation (2.3)).

**Lemma 4.1** *Let $F : \Sigma \times [0, T) \mapsto M$ be a mean curvature flow in $M$ and $S$ be a parallel two-tensor on $M$. $\nabla$ denotes the connection on $M$. For any tangent vector of $M$, $(\cdot)^T$ denotes the tangential part in $T\Sigma$ and $(\cdot)^\perp$ is the normal part in $N\Sigma$. Then*

$$(\frac{d}{dt} - \Delta)S(X, Y) = S((\nabla_X H)^T, Y) + S(X, (\nabla_Y H)^T)$$
$$- S((\nabla_{e_k}(\nabla_{e_k} X)^\perp)^T, Y) - S(X, (\nabla_{e_k}(\nabla_{e_k} Y)^\perp)^T) \quad (4.1)$$
$$- 2S((\nabla_{e_k} X)^\perp, (\nabla_{e_k} Y)^\perp)$$

*for any $X, Y \in T\Sigma$ and any orthonormal basis $\{e_k\}$ for $T\Sigma$, where $\Delta$ is the rough Laplacian on two-tensors over $\Sigma$.*

Now back to our setting when $\Sigma$ is Lagrangian in $T^{2n}$. $\{Je_k\}$ forms an orthonormal basis for $N_p\Sigma$. We define the second fundamental form by

$$h_{kij} = \langle \nabla_{e_k} e_i, J(e_j) \rangle$$

Thus $(\nabla_{e_k} e_i)^\perp = h_{kij} J(e_j)$ and $(\nabla_{e_k} J(e_i))^T = -h_{ikl} e_l$. Denote

$$H_j = \langle H, J(e_j) \rangle$$

Thus $(\nabla_{e_i} H)^T = -H_l h_{ijl} e_j$.

Plug these into equation (4.1), we derive

$$(\frac{d}{dt} - \Delta)S(e_i, e_j) = -H_p h_{ipl} S(e_l, e_j) - H_p h_{jpl} S(e_i, e_l)$$
$$+ h_{pki} h_{pkl} S(e_l, e_j) + h_{pkj} h_{pkl} S(e_i, e_l)$$
$$- 2 h_{kil} h_{kjm} S(J(e_l), J(e_m))$$

Recall $S(X,Y) = \langle J\pi_1(X), \pi_2(Y)\rangle$, so $S(J(e_l), J(e_m)) = \langle J\pi_1(J(e_l)), \pi_2(J(e_m))\rangle$. Since $J\pi_1 = \pi_2 J$ and $J\pi_2 = \pi_1 J$, we derive

$$S(J(e_l), J(e_m)) = \langle JJ\pi_2(e_l), J\pi_1(e_m)\rangle = -\langle \pi_2(e_l), J\pi_1(e_m)\rangle = -S(e_m, e_l) = -S(e_l, e_m)$$

The last step is because $S(\cdot, \cdot)$ is symmetric on any Lagrangian subspace.

Therefore, we obtain

$$(\frac{d}{dt} - \Delta)S_{ij} = (h_{pki}h_{pkl} - H_p h_{ipl})S_{lj} + (h_{pkj}h_{pkl} - H_p h_{jpl})S_{il} + 2h_{kil}h_{kjm}S_{lm} \tag{4.2}$$

Now $h_{pki}h_{pkl} - H_p h_{ipl} = R_{il}$ is indeed the Ricci curvature on $\Sigma$. This equation is also derived in [14]. Since $h_{kil}h_{kjm}S_{lm}$ is positive definite if $S_{ij}$ is, the positivity of $S_{ij}$ being preserved is a direct consequence of Hamilton's maximum principle for tensors [3]. To obtain the $C^2$ bound of $u$, we recall the following Lemma from [14]

**Lemma.** *Given any $\epsilon > 0$, the condition $S_{ij} - \epsilon g_{ij} > 0$ is preserved along the mean curvature flow.*

To see what this means in terms of the eigenvalues of $D^2 u$, we choose a particular orthonormal basis for $T_p\Sigma$ at a point $p$ that we are interested. The tangent space of $\Sigma$ is the graph of $df : \mathbb{R}^n \mapsto \mathbb{R}^n$. Recall the complex structure $J$ is chosen so that the target $\mathbb{R}^n$ is the image under $J$ of the domain $\mathbb{R}^n$. Because $\langle df(\cdot), J(\cdot)\rangle$ is symmetric, we can find an orthonormal basis $\{a_i\}_{i=1\cdots n}$ for the domain $\mathbb{R}^n$ so that

$$df(a_i) = \lambda_i J(a_i)$$

Then

$$\{e_i = \frac{1}{\sqrt{1+\lambda_i^2}}(a_i + \lambda_i J(a_i))\}_{i=1,\cdots,n} \tag{4.3}$$

becomes an orthonormal basis for $T_p\Sigma$ and $\{J(e_i)\}_{i=1\cdots n}$ an orthonormal basis for the normal bundle $N_p\Sigma$.

We compute for each $i$,

$$S(e_i, e_i) = \langle J\pi_1(e_i), \pi_2(e_i)\rangle = \frac{\lambda_i}{1+\lambda_i^2}$$

Now $\frac{\lambda_i}{1+\lambda_i^2} > \epsilon$ implies a uniform upper bound on $\lambda_i's$, the eigenvalues of $D^2u$.

An alternative way to get the $C^2$ bound is to consider the Gauss map of the mean curvature flow. By definition 2.2, $D^2u$ is represented by $\langle df(a_i), J(a_j) \rangle$. Comparing with equation (3.1), $z_{ij} = \langle df(a_i), J(a_j) \rangle$, therefore

$$\phi_0 \circ \gamma = -\frac{1}{2} \ln \det(I + (D^2u)^2)$$

where $\phi_0 = -\frac{1}{2} \ln \prod(1+\lambda_i^2)$ is defined in Proposition 3.2 and $\gamma : \Sigma \mapsto LG(n)$ is the Gauss map.

By equation (3.8) and Theorem A in [23], $-\ln \det(I + (D^2u)^2)$ is a supersolution of the nonlinear heat equation under the condition $\lambda_i > 0$, or

$$(\frac{d}{dt} - \Delta) \ln \det(I + (D^2u)^2) \leq 0$$

By the maximum principle, $\sup_t \ln \det(I + (D^2u)^2)$ is non-increasing in $t$. Since $\ln \det(I + (D^2u)^2) = \frac{1}{2} \ln \prod(1 + \lambda_i^2)$, all $\lambda_i$ are uniformly bounded. From here we conclude a uniform $C^2$ bound of $u$.

To prove the long time existence, recall the explicit formula from [20] (see also [17] for the elliptic version) for $*\Omega = \frac{1}{\sqrt{\prod(1+\lambda_i^2)}}$.

$$(\frac{d}{dt} - \Delta)(\ln *\Omega) = \left\{ \sum_{i,j,k} h_{ijk}^2 + \sum_{k,i} \lambda_i^2 h_{iik}^2 + 2 \sum_{k,i<j} \lambda_i \lambda_j h_{ijk}^2 \right\} \quad (4.4)$$

Also $h_{ijk} = \langle \nabla_{e_i} e_j, J(e_k) \rangle$ denotes the second fundamental form with respect to the basis in equation (4.3). The original formula in [20] is for $*\Omega$, to convert into the present form, recall $(*\Omega)_k = - *\Omega(\sum_i \lambda_i h_{iik})$ from [17].

We can now prove the long time existence as in [20]. By the positivity of $\lambda_i$ and equation (4.4), we obtain

$$(\frac{d}{dt} - \Delta)(\ln *\Omega) \geq |A|^2 \quad (4.5)$$

We integrate this inequality against the backward heat kernel and study the blow-up behavior at any possible singular points. A crucial point is any function defined on the Grassmannian is invariant under blow-up. The right hand side $|A|^2$ helps us to conclude any parabolic blow-up limit is totally geodesic and long time existence follows from White's regularity theorem [25].

From equation (2.2) we see that $\frac{du}{dt} = \alpha$ is given by the Lagrangian angle. On the other hand it is well known [12] that for the parametric mean curvature flow the Lagrangian angle satisfies the evolution equation

$$\frac{d}{dt}\alpha = \Delta\alpha \tag{4.6}$$

so that the maximum principle implies a uniform bound of $\frac{du}{dt}$. This means that $u$ is bounded in $C^1$ with respect to the time variable. We can even prove uniform $C^2$ bounds in time since by Theorem 1.3 in [14] $|H|$ and $|d^\dagger H| = |\Delta\alpha|$ are uniformly bounded.

An alternative proof of the long time existence which also implies convergence is to utilize the $C^{2,\alpha}$ estimate for nonlinear parabolic equations by Krylov [6] or [5](see section 5.5). To apply it, we still need to check the concavity of

$$\frac{1}{\sqrt{-1}} \ln \frac{\det(I + \sqrt{-1}D^2u)}{\sqrt{\det(I + (D^2u)^2)}}$$

in the space of symmetric matrices with the flat metric. This can be checked using a Lemma of Caffarelli, Nirenberg and Spruck in [1] (section 3 page 276), see also [14].

With the $C^{2,\alpha}$ bound in space and the $C^{1,\alpha}$ bound in time, the convergence now follows from standard Schauder estimates and Simon's theorem [10]. Equation (4.5) then implies the limit is a flat Lagrangian submanifold.

## 5 Other equivalent conditions

As was remarked in [21] (section 2) and [17] (see the remark at the end of the paper), the condition $u$ being convex corresponds to a region $V$ on the Lagrangian Grassmannian. Since the geometry of a Lagrangian submanifold is invariant under the unitary group $U(n)$, Theorem A applies whenever the Gauss map of a Lagrangian submanifold lies in a $U(n)$ orbit of $V$. To be more precise, consider $S$ as a bilinear form defined on $\mathbb{C}^n \cong T(T^{2n})$, given any $U \in U(n)$ we may consider $S_U$ defined by

$$S_U(\cdot, \cdot) = S(U(\cdot), U(\cdot))$$

Notice that $JU = UJ$ as linear transformations on $\mathbb{C}^n$. It is not hard to see $S_U$ again defines a symmetric bilinear form on any Lagrangian subspace. Also

$$S(J(X), J(Y)) = -S(X, Y)$$

for any $X, Y$ in a Lagrangian subspace.

Now $F^*S_U > 0$ for $F : \Sigma \mapsto T^{2n}$ implies the submanifold $\Sigma$ can be locally written as a graph over a different Lagrangian plane with a convex potential function. The new Lagrangian plane is indeed the image of the domain $\mathbb{R}^n$ under $U$. This corresponds to choosing a different base point in the Lagrangian Grassmannian in §3.

**Corollary A** *Let $F : \Sigma \mapsto T^{2n}$ be a Lagrangian submanifold. Suppose there exists an $U \in U(n)$ such that $F^*S_U$ is positive definite on $\Sigma$. Then the mean curvature flow of $\Sigma$ exists for all time and converges smoothly to a flat Lagrangian submanifold.*

Suppose $\Sigma$ is the graph of $f : T^n \mapsto T^n$ then the condition $F^*S_U > 0$ can be expressed in terms of the eigenvalues of the potential function $u$.

Recall from [21], given any splitting of $\mathbb{C}^n$, an element $U \in U(n)$ can be represented by a $2n \times 2n$ block matrix

$$\begin{bmatrix} P & -Q \\ Q & P \end{bmatrix}$$

with

$$PP^T + QQ^T = I, -PQ^T + QP^T = 0$$

Now corresponds to the slitting $T(T^{2n}) = \mathbb{C}^n = T\Sigma \oplus N\Sigma$ and the bases (equation (4.3))$\{e_i = \frac{1}{\sqrt{1+\lambda_i^2}}(a_i + \lambda_i J(a_i)), Je_i\}_{i=1,\cdots,n}$, we have $Ue_i = \sum_k P_{ki}e_k + \sum_l Q_{li}Je_l$. Then

$$S_U(e_i, e_i) = \sum_k (P_{ki}^2 - Q_{ki}^2)\frac{\lambda_k}{1+\lambda_k^2} + \sum_k P_{ki}Q_{ki}(\frac{1-\lambda_k^2}{1+\lambda_k^2})$$

Therefore the positive definiteness of $S_U$ is the same as requiring the above expression to be positive for each $i$. Take $P = Q = \frac{1}{\sqrt{2}}I$ which amounts to rotating each complex plane by $\frac{\pi}{4}$, we obtain $S_U(e_i, e_i) = \frac{1}{2}(\frac{1-\lambda_i^2}{1+\lambda_i^2})$. Therefore we have **Corollary B** *Let $\Sigma$ be a Lagrangian submanifold in $T^{2n}$. Suppose $\Sigma$ is the graph of $f : T^n \mapsto T^n$ and the absolute values of the eigenvalues of the potential function $u$ are less than one. Then the mean curvature flow of $\Sigma$ exists for all time and converges smoothly to a flat Lagrangian submanifold. In particular, during the evolution the absolute values of all eigenvalues stay less than one.*

That the flow preserves the property of $u$ having eigenvalues of absolute value less than one was also shown in [13], Theorem 2.6.3.